\documentclass[12pt]{article}
\usepackage{ae,amsfonts,euscript,enumerate}

\newcommand{\bA} {\mathbb{A}}
\newcommand{\bZ} {\mathbb{Z}_p}

\newcommand{\bQ}{\mathbb{Q}_p}

\newcommand{\ZZ}{\mathbb{Z}}

\newcommand{\bs}{{\bf s}}
\newcommand{\bx}{{\bf x}}

\newcommand{\sRS}{{\cal R}}
\newcommand{\sA}{{\cal A}}

\newtheorem{thm}{Theorem}[section]
\newtheorem{exam}{Example}
\newtheorem{defn}{Definition}[section]

\newtheorem{prop}{Proposition}[section]
\newtheorem{lem}{Lemma}[section]
\newtheorem{cor}{Corollary}[section]

\newtheorem{quest}{Question}[section]

\def\eqref#1{(\ref{#1})}
\def\qed{~\vrule height8pt width 5pt depth -1pt\medskip}

\begin{document}
\bibliographystyle{plain}
\begin{center}
{\bf \large A generalised Skolem-Mahler-Lech theorem for affine varieties
}\\
\vspace{7 mm} \centerline{Jason P. Bell}

\centerline{Department of Mathematics}
\centerline{
Simon Fraser
University}
\centerline{
8888 University Dr.}
\centerline{Burnaby, BC V5A 1S6.}
\centerline{ CANADA}
\centerline{\tt jpb@math.sfu.ca}
\end{center}
\begin{abstract} The Skolem-Mahler-Lech theorem states that
if $f(n)$ is a sequence given by a linear recurrence over a field of
characteristic $0$,
then the set of $m$ such that $f(m)$ is equal to $0$ is the union of
a finite number 
of arithmetic progressions in $m\ge 0$ and a finite set.
We prove that if $X$ is a subvariety of an affine variety $Y$ over a field
of characteristic $0$ and ${\bf q}$ is a
 point in
$Y$,
and $\sigma$ is an automorphism of $Y$, then
the set of $m$ such that $\sigma^m({\bf q})$ lies in $X$ is a union of
a finite number 
of complete doubly-infinite arithmetic progressions and a finite
set.  We show that this is a generalisation of the Skolem-Mahler-Lech
theorem.
\end{abstract}
key words: automorphisms, affine space, Skolem-Mahler-Lech theorem \\
Primary classification: 11D45 \\
Secondary classification: 14R10, 11Y55, 11D88  
\section{Introduction}
The Skolem-Mahler-Lech theorem is a well-known result 
which describes the set of solutions in $n$ to the equation $f(n)=0$, where
$f(n)$ is a sequence given by a linear recurrence.
\begin{thm} (Skolem-Mahler-Lech)
Let $f(n)=\sum_{j=1}^r a_jf(n-j)$ for $n\ge r$ and $a_1,\ldots ,a_r$ lie
in a field of characteristic $0$.  
Then the set of $m$ such that $f(m)=0$ is a union of a finite number
of arithmetic \label{thm: SMLlr}
progressions in $m\ge 0$ and a finite set. 
\end{thm}
This theorem was first proved for linear recurrences over the rational
numbers by Skolem \cite{Sk}.  It was next proved for linear recurrences 
over the algebraic numbers by Mahler \cite{Mah1}.
The version
above was proven first by Lech \cite{Lech} and later 
by Mahler \cite{Mahler}, \cite{Maherror}.
This history of this theorem can be found in the book by 
Everest et al. \cite{EV}.
\vskip 2mm
There are many different proofs and extensions of 
the Skolem-Mahler-Lech theorem in the
literature \cite{B},\cite{Han},\cite{EV},\cite{vT},\cite{vdP}.  
These proofs all use $p$-adic methods in
some way, although the result is valid in any field of characteristic $0$.

In Section \ref{sec: 3}, we show that the Skolem-Mahler-Lech theorem can
be rephrased in the following linear algebraic form.       
\begin{thm} \label{thm: SMLlin} Let $K$ be a field of characteristic $0$, and let
 $\sigma:K^n\rightarrow K^n$ be an invertible
linear map,
let ${\bf v}\in K^n$, and let $W$ be a vector subspace of
$K^n$ of codimension $1$.  Then the set of $m\in {\mathbb N}$ such that
 $\sigma^m({\bf v})\in W$ is the union of a finite number of arithmetic
progressions in $m\ge 0$ and a finite set.
\end{thm}

A related version of the Skolem-Mahler-Lech theorem was given by
Denis
\cite[Th\'eor\`eme 1]{Denis}.
We give an ``algebro-geometric'' generalisation of Theorem \ref{thm: SMLlin}
which we
now describe.
Given a ring $R$ with unit, we let $\bA_R^n$ denote $R^n$.
We consider the \emph{polynomial
automorphisms}
of $\mathbb{A}^n_R$, ${\rm Aut}(\bA_R^n)$, which are defined as follows.
\begin{defn} Let $R$ be a commutative ring with unit.  We define 
${\rm Aut}(\bA_R^n)$ to be the set of invertible polynomial maps
$\sigma=(f_1,\ldots ,f_n)\in R[x_1,\ldots ,x_n]^n$ with the property that
$\sigma^{-1}=(g_1,\ldots ,g_n)$ for some $g_1,\ldots ,g_n\in R[x_1,\ldots ,x_n]$.
An element $\sigma\in {\rm Aut}(\bA_R^n)$ is called an \emph{automorphism}
of $\bA_R^n$.
\end{defn}
When $R$ is
a field,
linear automorphisms are a special case of this in which we insist that
all polynomials involved in defining the maps are homogeneous of degree $1$;
i.e., they are invertible linear maps.

Our main result is the following algebro-geometric generalisation of this
result, to affine varieties $Y$.
\begin{thm} Let $Y$ be an affine variety over a field $K$ of characteristic $0$, let
${\bf q}$ be a point in $Y$ and let $\sigma$ be an automorphism of
$Y$. 
If $X$ is a subvariety of $Y$,
then the set \label{thm: intro}
$\{ m\in {\mathbb Z}~|~
\sigma^{m}({\bf q})\in X\}$ is a union of a finite number of complete
doubly-infinite arithmetic
progressions and a finite set.
\end{thm}

Theorem \ref{thm: SMLlin} is a special case of 
Theorem \ref{thm: intro} in which $Y=\bA_K^n$, $X$ is a subspace
of $Y$ of codimension $1$ and $\sigma$ is a linear automorphisms.  
While linear automorphisms are well-understood, the
polynomial
automorphism group of $\mathbb{A}^n_K$ remains a mystery in general (see
van den Essen \cite{Ess}).  Polynomial
automorphisms 
form a very large set.  For $n=2$, there is a
concrete description of ${\rm Aut}(\mathbb{A}_K^2)$
given by van der Kulk \cite{Kulk}.
There are several proofs of van der Kulk's result, see for example,
Gizatullin and Danilov \cite{GD}.
There are many open questions concerning the automorphisms of $\bA_K^n$.
Arguably, the most famous of these questions is the Jacobian conjecture
\cite{Ess}.  Whether Theorem \ref{thm: intro} has any application to the
Jacobian conjecture remains to be seen.  

We show elsewhere that Theorem \ref{thm: intro} can be further 
generalised to certain non-affine algebraic varieties, 
including $\mathbb{P}^n_K$, Fano varieties, and certain projective surfaces 
\cite{Be}.  In this context the problem is related to whether certain subrings
of twisted homogeneous coordinate rings of varieties are Noetherian \cite{Ro}.

The main part of the
proof of Theorem \ref{thm: intro} is for the case $Y=\bA^n_K$, given as
Theorem \ref{thm: introA} in \S 4.  The general case follows from the case
$Y=\bA^n_K$, using a theorem of Srinivas \cite{Sr}, and is given in \S 5.
The proof of Theorem \ref{thm: introA} proceeds by a reduction to the
$p$-adic case, given in \S 3.  This 
is the same reduction
 used in other proofs of the Skolem-Mahler-Lech
 theorem.  We then show that for a point ${\bf q}$ of
$\bA_{\bZ}^n$ and $\sigma\in {\rm Aut}(\bA_{\bZ}^n)$, there exist $p$-adic analytic
maps $f_1,\ldots ,f_n$ and some $j\ge 1$ such that for some $i<j$,
$$\sigma^{jm+i}({\bf q})=(f_1(m),\ldots ,f_n(m))$$ for all $m\in {\mathbb Z}$.
We think of this as embedding the orbit of ${\bf q}$ under $\sigma$ in a $p$-adic
analytic arc.  This is the key idea in the proof.
Having done this, we complete the proof of Theorem \ref{thm: introA} by using
 Strassman's theorem \cite{St} (see also Cassels \cite[Theorem 4.1, p. 62]{C}), which
says that a $p$-adic analytic map from $\mathbb{Z}_p$ to itself having infinitely
many zeros in $\mathbb{Z}_p$ must be identically zero.

\section{The linear algebraic form of the Skolem-Mahler-Lech theorem} \label{sec: 3}
The Skolem-Mahler-Lech theorem can also be stated as a theorem about coefficients of rational power
series over a field of characteristic $0$.
\begin{thm} (Skolem-Mahler-Lech, second version)
Let $\sum_{m=0}^{\infty} c_m z^m$ be the power series expansion of a
rational function over a field of characteristic $0$.
Then the set of $m$ such that $c_m=0$ is a union of a finite number
of arithmetic \label{thm: SML}
progressions in $m\ge 0$ and a finite set.
\end{thm}
This is the form given by Mahler \cite{Mah1},\cite{Mahler}.
We proceed to show that this version of the Skolem-Mahler-Lech theorem 
is equivalent to a statement in terms of linear algebra,
given as Theorem \ref{thm: SMLlin}.
This equivalence arises from the following theorem of
Sch\"utzenberger
\cite{Sc} (see, also, 
Hansel \cite[Proposition 1.1]{Han}).
\begin{thm} (Sch\" utzenberger)
 Let $f(x)\in K[[x]]$ be a rational function.  Then there exists
an invertible matrix $M$, column vectors ${\bf v},{\bf w}$, and a polynomial
$P(x)$ such that
$$f(x) = P(x) + \sum_{i=0}^{\infty} \big({\bf v}^{\rm T}M^i{\bf w}\big)x^i.$$
Conversely, given an invertible matrix $M$ and column vectors ${\bf v},{\bf w}$,
the power series
$$g(x) \ = \ \sum_{i=0}^{\infty} \big({\bf v}^{\rm T}M^i {\bf w}\big)x^i$$ is rational.
\label{thm: equivalent}
\end{thm}
We now show Theorem \ref{thm: SML} and \ref{thm: SMLlin} are easily deduced
from each other.  In this sense we regard them as equivalent.
\vskip 2mm
\noindent
{\bf Proof that Theorem \ref{thm: SML} $\iff$ Theorem \ref{thm: SMLlin}.}
Let $K$ be a field of characteristic $0$.
By fixing a basis for $K^n$, we may regard a linear automorphism $\sigma$
of $K^n$ as being given by multiplication by an invertible matrix $M$.
Furthermore, a codimension $1$ vector subspace $W$ of $K^n$ is of the form
$$\{{\bf x}~|~{\bf v}^{\rm T}
{\bf x}=0\}$$ for some vector ${\bf v}$.  Hence the set of all $m\in {\mathbb
N}$ such that $\sigma^m({\bf w})\in W$ is exactly the set of
$m\in \mathbb{N}$ such that ${\bf v}^{\rm T}M^m{\bf w}=0$.  By Theorems
\ref{thm: equivalent} and
\ref{thm: SML}, we deduce Theorem \ref{thm: SMLlin}.  

Conversely, if we have
a rational power series $f(x)$, then by Theorem \ref{thm: equivalent},
there is an invertible
matrix $M$ and column vectors ${\bf v}$ and ${\bf w}$ such that
$$f(x)=P(x)+\sum_{i=0}^{\infty} \big({\bf v}^{\rm T}M^i{\bf w}\big)x^i,$$ 
for some polynomial
$P(x)$.  Thus from Theorem \ref{thm: SMLlin} we immediately deduce the
Skolem-Mahler-Lech theorem.  \qed
\section{Embedding in ${\mathbb Z}_p$}
\label{sec: p-adic}
The results of this section are preparation to proving Theorem \ref{thm: intro} for 
$Y=\bA^n_K$ in \S 4.
It reduces the problem to working over a $p$-adic field using
an embedding argument due to Lech \cite{Lech}, which can be regarded as a
$p$-adic analogue of the Lefschetz principle.
\begin{lem} Let $K$
be a finitely generated
extension of $\mathbb{Q}$ 
and
let $\mathcal{S}$ be a finite subset of
$K$.  Then there is a prime $p$ such that
$K$ embeds in $\mathbb{Q}_p$; moreover, we can choose this embedding so
that every element of $\mathcal{S}$ is sent to an element in $\mathbb{Z}_p$.
\label{lem: p} \end{lem}
{\bf Proof.}  We assume, without loss of generality, that
 $0\not\in \mathcal{S}$.
Write
$K=
\mathbb{Q}(t_1,\ldots ,t_d)(\theta)$ , where
$t_1,\ldots ,t_d$ are algebraically
independent over $\mathbb{Q}$ and $\theta$ is algebraic over
$\mathbb{Q}(t_1,\ldots ,t_d)$.  We note this can be done by the primitive
element theorem.
Let $$f(x) \ = \ f(t_1,\ldots ,t_{d};x) \ \in \
 {\mathbb Q}(t_1,\ldots ,t_{d})[x]$$ be the minimal polynomial
of $\theta$.  By clearing denominators if necessary, we may assume that
$f(x)\in {\mathbb Z}[t_1,\ldots ,t_{d}][x]$.
Let $\Delta(t_1,\ldots ,t_{d})
\in {\mathbb Z}[t_1,\ldots ,t_{d}] $ denote the
discriminant of $f(x)$ with respect to the variable $x$.  For each
$s\in \mathcal{S}$ there is some 
 polynomial
$g_{s}
(x)\in {\mathbb Q}(t_1,\ldots ,t_{d})[x]$ satisfying $g_{s}(\theta)=s$.
Choose a polynomial $B_{s}(t_1,\ldots ,t_{d})\in {\mathbb Z}[t_1,\ldots ,t_d]$
such that
$B_sg_s\in {\mathbb Z}[t_1,\ldots ,t_{d}][x]$.   Let
$A_s\in {\mathbb Z}[t_1,\ldots ,t_d]$ denote the resultant of $f(x)$ and
$B_sg_s(x)$ with respect to the variable $x$.
  Choose integers $a_1,\ldots ,a_{d}$ such that:
\begin{enumerate}
\item{$\Delta(a_1,\ldots ,a_{d})\not =0$;}
\item{$f(a_1,\ldots ,a_d,x)$ is a non-constant polynomial;}
\item{
$  B_{s}(a_1,\ldots ,a_{d})\not = 0$ for $s\in
 \mathcal{S}$; and}
\item{$A_s(a_1,\ldots ,a_d)
\not = 0$
for all $s\in \mathcal{S}$.}
\end{enumerate}
Choose some prime $p$ which satisfies:
\begin{enumerate}
\item{$B_{s}(a_1,\ldots ,a_{d})\not \equiv 0 ~(\bmod~ p)$ for all
$s\in \mathcal{S}$;}
\item{$\Delta(a_1,\ldots ,a_{d})\not \equiv 0 ~(\bmod~ p)$;}
\item{$A_s(a_1,\ldots ,a_d)\not \equiv 0~(\bmod~ p)$; and}
\item{$f(a_1,\ldots ,a_{d};x)$ has a root mod $p$.}
\end{enumerate}
The Chebotarev
density theorem \cite[Theorem 10, p. 109]{La}
 shows that such a prime must exist.
Pick algebraically independent elements $\mu_1,\ldots ,\mu_{d}$ of
${\mathbb Z}_p$.  This can be done since $\mathbb{Z}_p$ is uncountable.  
Notice that the function
$f(a_1+p\mu_1,\ldots ,a_{d}+p\mu_{d},x)$ has a root mod $p$.  Since
$\Delta\not =0~(\bmod~ p)$, we see by Hensel's lemma
that there is a number $\tilde{\theta}$ in
${\mathbb Z}_p$ such that
$$f(a_1+p\mu_1,\ldots ,a_{d}+p\mu_{d},\tilde{\theta}) \ = \ 0.$$
We therefore have an embedding
$K \hookrightarrow {\mathbb Q}_p$ given by
$$t_i\mapsto a_i + p\mu_i,\qquad \theta \mapsto \tilde{\theta}.$$
Moreover, the conditions on $A_s$ and $B_s$ for $s\in \mathcal{S}$ ensure
that elements of $\mathcal{S}$ get sent to
elements of ${\mathbb Z}_p$ in this embedding.
\qed
\vskip 2mm
We now describe how Lemma \ref{lem: p} will be applied to our situation.
We are working with an automorphism of $\mathbb{A}^n_K$ over a field
of characteristic $0$.  We are looking at the orbit of a point ${\bf q}=(q_1,
\ldots ,q_n)$ under
$\sigma$ and asking where this orbit intersects some subvariety $X$.  
Write $\sigma=(H_1,\ldots ,H_n)$ and $\sigma^{-1}=(G_1,\ldots,G_n)$
 for some polynomials 
$H_1,\ldots ,H_n,G_1,\ldots, G_n \in K[x_1,\ldots ,x_n]$ and write
${\bf q}=(q_1,\ldots ,q_n)$.  
Finally, choose polynomials $F_1,\ldots ,F_m$
such that $X$ is the zero set of the polynomials $F_1,\ldots ,F_m$.
Let $\mathcal{S}$ be the
set consisting of $q_1,\ldots ,q_n$ along with the coefficients
$H_1,\ldots ,H_n,G_1,\ldots ,G_n, F_1,\ldots ,F_m$.  It is no loss of 
generality to assume that $K$ is the finitely generated extension of 
$\mathbb{Q}$ generated by $\mathcal{S}$.  It follows from Lemma \ref{lem: p}
 that
there is some prime $p$ and an embedding $K\hookrightarrow {\mathbb Q}_p$
in which the elements of $\mathcal{S}$ get sent to elements in $\mathbb{Z}_p$.
Thus we may assume that $K={\mathbb Q}_p$ and that
$\sigma$ restricts to an automorphism of
$\bA_{\bZ}^n$; that is, $\sigma$ is bijective map 
from $\bA_{\bZ}^n$ to itself in which $\sigma$ and $\sigma^{-1}$ are both given
by $n$-tuples of polynomials in $n$ variables with the coefficients 
in $\bZ$.  Note that $\sigma$ necessarily has Jacobian matrix whose determinant is
a constant of $p$-adic absolute value $1$.
\section{Proofs for affine space}
In this section we prove the following theorem.
\begin{thm} Let $K$ be a field of characteristic $0$, let
${\bf q}$ be a point in $\bA_K^n$ and let $\sigma$ be an automorphism of
$\bA_K^n$.
If $X$ is a subvariety of $\bA^n_K$,
then the set \label{thm: introA}
$\{ m\in {\mathbb Z}~|~
\sigma^{m}({\bf q})\in X\}$ is a union of a finite number of complete
doubly-infinite arithmetic
progressions and a finite set.
\end{thm}

In light of section \ref{sec: p-adic}, 
we may assume that
we have an automorphism $\sigma$ of $\bA_{\bZ}^n$ and that $X$ is
a subvariety of $\bA_{\bZ}^n$ and that ${\bf q}$
is a point in $\bA_{\bZ}^n$.
Up to this point there have been no substantive differences between our proof
and the proof of the Skolem-Mahler-Lech theorem given by Lech \cite{Lech}.  At this
point, however, our arguments differ.  Lech used the fact that the
coefficients of a rational function are linear combinations of terms of the form
$\alpha^n n^j$, with $j\in \mathbb{N}$ and $\alpha\in \mathbb{C}$.  It is not
difficult to embed the function $\alpha^n n^j$ in an analytic arc in $\mathbb{C}$;
for instance, the function $f(z)=z^j \alpha^z$ is an entire function with $f(n)=\alpha^n n^j$.  In the $p$-adic case, one must be more careful to ensure that 
$f(z)$ is analytic on $\bZ$.  It is here that arguments involving arithmetic 
progressions come into play.  For general 
automorphisms $\sigma$ of $\bA^n_{\bZ}$, there is no nice
formula giving $\sigma^n({\bf q})$ and we must do more work to embed 
$\{\sigma^n({\bf q})\}$ in a 
$p$-adic analytic 
arc.  This is accomplished via the next two lemmas.   

\begin{lem} Let 
$\sigma=(F_1,\ldots ,F_n):\bA_{\bZ}^n\rightarrow \bA_{\bZ}^n$ be a surjective 
polynomial map 
with $F_1,\ldots ,F_n\in \bZ[x_1,\ldots ,x_n]$ and which has the property that
its Jacobian determinant is a constant and has $p$-adic norm $1$.
Then there is a positive integer $j$ such that
$\sigma^j=(H_1,\ldots ,H_n)$ has the following two properties:
\begin{enumerate}
\item{ $\sigma^j(s_1,\ldots ,s_n)=(s_1',\ldots ,s_n')$ satisfies
$s_i\equiv s_i'~(\bmod~ p)$ for $1\le i\le n$ for all $(s_1,\ldots ,s_n)\in 
\bZ^n$;}
\item{
The Jacobian of $(H_1,\ldots ,H_n)$ evaluated at
$x_1=s_1,\ldots ,x_n=s_n$ is of the form $I+pM$, for some matrix $M$ with
entries in $\bZ$.}
\end{enumerate}
\label{lem: j}
\end{lem}
 {\bf Proof.}  
We have 
$\sigma(s_1,\ldots ,s_n) = (F_1(s_1,\ldots ,s_n),\ldots ,F_n(s_1,\ldots ,s_n))$ 
for some $F_1,\ldots ,F_n\in {\mathbb Z}_p[x_1,\ldots ,x_n]$.
Since the $F_i$'s have coefficients in ${\mathbb Z}_p$,
we see that $\sigma(s_1,\ldots ,s_n)
=(r_1,\ldots ,r_n)$ for some $r_1,\ldots ,r_n$
in ${\mathbb Z}_p$.  Moreover, if $$(r_1,\ldots ,r_n) \ = \ \sigma(s_1,
\ldots ,s_n)\qquad {\rm and}\qquad (r_1',\ldots ,r_n') \ = \ 
\sigma(s_1',\ldots ,s_n'),$$ and 
$s_i\equiv s_i'~(\bmod~ p\bZ)$ for $1\le i\le n$
then $r_i\equiv r_i'~(\bmod~ p\bZ)$ for $1\le i\le n$.
Hence there exist $k>\ell>0$ such
that
$\sigma^k(s_1,\ldots ,s_n)$ and
$\sigma^{\ell}(s_1,\ldots ,s_n)$ have equivalent coordinates mod $p$ for
each $(s_1,\ldots ,s_n)\in \bZ^n$.
Since $\sigma^{\ell}$ is surjective, we see that
$\tau=\sigma^{k-\ell}$ has the property that
if $\tau(s_1,\ldots ,s_n)=(s_1',\ldots ,s_n')$, then
$s_1'\equiv s_1,\ldots ,s_n'\equiv s_n~(\bmod~ p)$.  Hence any power of 
$\tau$ has this property.  Since $\sigma$ has a constant Jacobian determinant in
$\bZ^{\times}$, so
does $\tau$.
Write $J(\mu,{\bf s})$ for the Jacobian of a polynomial map
$\mu$
 at
a point ${\bf s}\in \bA_{\bZ}^n$.  
Let $m$ denote the order of ${\rm GL}_n(\mathbb{F}_p)$.  Then for any point ${\bf s}=(s_1,\ldots ,s_n)\in \bA_{\bZ}^n$ we
have 
$$J(\tau^m,{\bf s}) \ = \ J(\tau,\tau^{m-1}({\bf s}))J(\tau,\tau^{m-2}({\bf s}))\cdots
J(\tau,{\bf s}).$$
Note that $\tau^i({\bf s})$ has coordinates congruent to the coordinates of
${\bf s}$ mod $p\bZ$.  Since $J$ has entries given by polynomials with coefficients in $\bZ$, we conclude that $J(\tau,\tau^i({\bf s}))\equiv J(\tau,{\bf s})
~(\bmod ~p\bZ)$, where $\equiv$ is taken coordinate-wise.  Hence
$$J(\tau^m,{\bf s})\equiv J(\tau,{\bf s})^m \equiv I_n ~(\bmod~ p\bZ).$$
Thus $\tau^m$ satisfies the conditions in the statement of the lemma.
It follows that there is some $j$ such that $\sigma^j$
satisfies the conditions given in the statement of the lemma. \qed 
\vskip 1mm
\noindent
To avoid any possible confusion in interpreting the preceding lemma, 
we make the remark that if a polynomial $f(x_1,\ldots ,x_n)\in \bZ[x_1,\ldots
,x_n]$ has the property that $f(s_1,\ldots ,s_n)\equiv 0~(\bmod ~
p\bZ)$ for all $(s_1,\ldots, s_n)\in \bZ^n$, then it not necessarily the case
that the coefficients of $f$ are all in $p\bZ$.  For example, the polynomial 
$f(x)=x^p-x$ has this property in the case $n=1$.

\begin{lem} \label{le21}
Let $p\ge 2$ be  prime.
Let $N \ge 1$ be a natural number, and let
\[ S_N \ := \ \Bigg\{ c + \sum_{i=1}^N p^i h_i(z)~|~c\in \bZ, h_i(z)\in
\sRS, {\rm deg}(h_i)\le 2i-1\Bigg\}\]
and
\[ T_N \ := \ S_N+\Bigg\{ \sum_{i=1}^M p^i h_i(z)~|~M\ge 1,  h_i(z)\in
\sRS, {\rm deg}(h_i)\le 2i-2\Bigg\}.\]
Then the $\bZ$-subalgebra $\sA[S_N]$ of $\sRS$ generated by $S_N$ is contained in $T_N$.
\label{lem: SN}
\end{lem}

\noindent {\bf Proof.} Since $S_N$ and $T_N$ are both closed under addition and
$S_N\subseteq T_N$, it is sufficient to show that $S_NT_N\subseteq T_N$.
To do this, suppose
\[H(z) \ = \ c + \sum_{i=1}^N p^i h_i(z) \ \in \ S_N\] and
\[G(z) \ = \ d + \sum_{i=1}^M p^i g_i(z) \ \in \ T_N,\]
where $c,d\in \ZZ_p$ and $h_i(z),g_i(z)\in \sRS$ with ${\rm deg}(g_i)\le
2i-2$ for $i>N$ and ${\rm deg}(g_i),{\rm deg}(h_i)\le 2i-1$ for $i\le N$.
We must show that $H(z)G(z)\in T_N$.  Notice that
\[
H(z)G(z) \ = \ cG(z)+dH(z) -cd + \sum_{i=1}^N p^i h_i(z)\sum_{j=1}^M p^j
g_j(z). \]
Then $cG(z), dH(z)$, and $cd$ are all in $T_N$ and since $T_N$ is closed
under addition, it is sufficient to show that
\[\sum_{i=1}^N p^i h_i(z)\sum_{j=1}^M p^j g_j(z)  \ = \ \sum_{k=2}^{N+M}
p^k \sum_{i=1}^{k-1}g_i(z)h_{k-i}(z)\]
is in $T_N$.  But since $g_i(z)$ has degree at most $2i-1$ and
$h_{k-i}(z)$ has degree at most $2(k-i)-1$, we see that 
$$\sum_{i=1}^{k-1} g_i(z)h_{k-i}(z)$$ has degree at most $2k-2$.  It
follows that
\[\sum_{i=1}^N p^i h_i(z)\sum_{j=1}^M p^j g_j(z)  \ \in \ T_N\]
as Mahler series.
The result follows. \qed

%
%

\begin{thm} ($p$-adic Analytic Arc Theorem) \label{lem: 2nd} 
Let $p\ge 5$ be prime, and let
$\bs =(s_1,\ldots ,s_d)\in \mathbb{Z}_p^n$.
Suppose $\tau=(H_1, H_2, ..., H_n)$ is a polynomial map
with each $H_i(x_1, ..., x_n) \in{\mathbb Z}_p[x_1,\ldots ,x_n]$,
which satisfies:\\

(i) $H_i(s_1,\ldots ,s_n)\equiv s_i~(\bmod ~ p\mathbb{Z}_p)$ for $1\le i\le n$;\\

(ii) The Jacobian matrix  $J(\tau; \bx)= [ \frac{\partial H_i}{\partial x_j}]$ 
evaluated at $\bx=\bs\in {\mathbb Z_p}^n$ is of the form
$I_n + pM$, i.e. it is congruent to the identity $(\bmod ~p\bZ)$.\\

Then there exist power series 
$f_1(z),\ldots ,f_n(z) \in \bQ[[z]]$ 
which are analytic on $\ZZ_p$  and
which satisfy: \\

(1) $f_i(z+1) \ = \
H_i(f_1(z),\ldots ,f_n(z))$ for $z\in \mathbb{Z}_p$ for $1\le i\le n$; \\

(2) $f_i(0)= s_i$ for $1\le i\le n$.

\label{lem: comp}
\end{thm}

\noindent {\bf Proof.}
We will construct $(f_1(z), ..., f_{n}(z))$ by approximation
$(\bmod~ p^{j} \sRS)$.  To initialize we define $g_{i,0}(z)=s_i,$  a
constant function,
for $0\le i\le n$. We will recursively define polynomials 
$h_{i,j}(z) \in \sRS$  ($1\le i\le d$)  for $j \ge 1$ and define
$$
g_{i,j}(z) :=s_i+\sum_{k=1}^{j} p^{k}h_{i,k}(z).
$$
We prove by induction on $j$ that for each $j \ge 1$ there
exist polynomials $h_{i,j}(z) \in \sRS$  ($1\le i\le n$) which satisfy
the three conditions:

(1) $h_{i,j}(0)=0$ for $1\le i\le n$;

(2) $h_{i,j}(z)$ has degree at most $2j-1$ for $j\ge 1$ and $1\le i\le n$;

(3) If $g_{i,j}(z)=s_i+\sum_{k=1}^{j} p^{k}h_{i,k}(z)$, then
  \[g_{i,j}(z+1)\equiv H_i(g_{1,j}(z),\ldots ,g_{n,j}(z)) ~(\bmod ~
p^{j+1}\sRS ).\]

For the base case of the induction
we take $j=0$ and set each $h_{i,0}(z)=0$, making the
convention that the zero polynomial has degree $-1$; then (1)-(3) hold.
 Now let $j\ge 1$ and assume that we have defined $h_{i,k}$ for $1\le i\le n$ 
and $0\le k <j $
so that conditions (1)-(3) hold.  
We use the notation of the statement of Lemma \ref{lem: SN}.  Then (2), (3) say that
$g_{i,k} \in S_{k}$ for $0 \le k \le j-1$.
Our object is to construct
\[
g_{i, j}(z) := g_{i,j-1}(z) + p^j h_{i, j}(z),
\]
in which polynomials $h_{i, j}(z) \in \sRS$ are to
be determined, so that conditions (1)-(3) hold.
By assumption
\begin{equation} \label{eq: gH}
g_{i,j-1}(z+1)-H_i(g_{1,j-1}(z),\ldots ,g_{d,j-1}(z))=p^jQ_{i,j}(z),
\end{equation}
with
$Q_{i,j}\in \sRS$
for $1\le i\le n$.
Conditions (2) and (3) show that
$g_{1,j-1}(z),\ldots ,g_{d,j-1}(z)$ are in $S_{j-1}$.
Thus we see that
\[ p^j Q_{ij}(z)\ = \ g_{i,j-1}(z+1)-H_i(g_{1,j-1}(z),\ldots ,g_{n,j-1}(z))\]
is in the $\ZZ_p$-algebra $\sA[S_{j-1}]$ generated by $S_{j-1}$ and hence 
by Lemma \ref{lem: SN} is in $T_{j-1}$.  It
follows that we can write
$p^j Q_{ij}(z)= c_{ij}+ \sum_{k=1}^M p^k q_{ijk}(z)$ for some $c_{ij}\in  
\ZZ_p$ and polynomials $q_{ijk}(z)\in \sRS$ such that ${\rm
deg}(q_{ijk})\le 2k-1$ for $k\le j-1$ and ${\rm deg}(q_{ijk})\le 2k-2$
for $k\ge j$.  Consequently,
$p^jQ_{ij}(z)$ is equivalent (mod $p^{j+1}\sRS$) to the polynomial
$$ c_{ij}+ \sum_{k=1}^{j} p^k q_{ijk}(z),$$
which is a polynomial of degree at most $2j-2$. 
Hence $Q_{ij}(z)$ is congruent (mod $p\sRS$) to a polynomial in $\sRS$ of
  degree at most $2j-2$.
To satisfy property
(3) for $j$ it  is sufficient to find 
$\{h_{i, j}(z) \in \sRS : 1 \le i \le n\}$ such
that 
\[
g_{i,j-1}(z+1)+p^j h_{i,j}(z+1)  - H_i(g_{1,j-1}(z)+p^jh_{1,j}(z),\ldots ,
g_{d,j-1}(z)+p^jh_{n,j}(z)) \in p^{j+1}\sRS
\]
for $1\le i\le d$.
Using the Taylor expansion and considering this expression modulo $p^{j+1}\sRS $, we see
\[
p^jQ_{i,j}(z) + p^j h_{i,j}(z+1)
  - p^j \sum_{\ell=1}^d h_{\ell,j}(z) \frac{\partial H_{i}}
{\partial x_{\ell}} (x_1, ..., x_n)
\Big|_{{x_{i}=g_{i,j}(z)}\atop{1 \le i \le n}}
\equiv 0 ~ (\bmod~p^{j+1} \sRS).
\]

It therefore suffices  to solve (mod~ $p \sRS$)
for $1 \le i \le n$  the system
\begin{equation} \label{eq: Q}
Q_{i,j}(z) + h_{i,j}(z+1)
  - \sum_{\ell=1}^d h_{\ell,j}(z) \frac{\partial H_{i}}
{\partial x_{\ell}}(x_1, ..., x_n) 
\Big|_{{x_{i}=g_{i,j}(z)}\atop{1 \le i \le n}}
  \equiv 0 ~(\bmod ~p\sRS),
\end{equation}
  where we may assume that $Q_{i,j}$ is of degree at most $2j-2$.
Now consider  the Jacobian matrix $M^{(j)}(z) \in M_{n \times n}(\sRS)$
with polynomial entries
\[
M^{(j)}(z)_{i\ell} := \frac{\partial H_{i}}
{\partial x_{\ell}}(x_1, ..., x_d) \Big|_{x_{1}=g_{1,j}(z),\ldots, x_n=
g_{n,j}(z)}
\]
Property (1) for $j$ yields for $1 \le i \le n$ that 
\[ 
g_{i, j}(z) \equiv s_i~ (\bmod~ p \sRS).
\]
It follows that
\[
M^{(j)}(z) \equiv J(\tau; \bs) ~(\bmod~ p \sRS).
\]
By hypothesis (ii) the matrix $M=  J(\tau; \bs) \in M_{n\times n}(\bZ)$
has $M\equiv I_n$ $(\bmod~p\bZ)$, where $I_n$
is the identity matrix, so we have
\[
M^{(j)}(z) \equiv M \equiv I_n     ~ (\bmod~ p \sRS).
\]
It follows that equation (\ref{eq: Q}) can be rewritten in the form 
\[ 
\left[\begin{array}{c} h_{1,j}(z+1) \\ \vdots \\ h_{n,j}(z+1)\end{array}
\right]
\equiv \left[\begin{array}{c} h_{1,j}(z) \\ \vdots \\ h_{n,j}(z)\end{array}
\right] - 
\left[\begin{array}{c} Q_{1,j}(z) \\ \vdots \\ Q_{n,j}(z)\end{array}
\right]
(\bmod~ p\sRS).
\]
Now 
there exists a  solution $[h_{1,j}(z), ..., h_{n, j}(z)] \in \sRS^d$ with 
$h_{i,j}(0)=0$ for $1\le i\le n$ and
$h_i(z)$ of degree at most $2j-1$, because each equation is separately
solvable in $\sRS$,
since
$$
\left( {{z+1}\atop{k+1}}\right) - \left( {{z+1}\atop{k}}\right) =\left(
{{z}\atop{k}}\right),~~~ k \ge 0.
$$ 
Thus conditions (1)-(3) are satisfied for $j$, completing the induction step.

We now set
\[
f_i(z) := s_i + \sum_{j=1}^{\infty} p^j h_{i, j}(z).
\]
Now each $h_{i,j}(z) \in \sRS$ is of degree at most $2j-1$ and hence
\[
h_{i,j}(z) = \sum_{k=0}^{2j-1}  c_{ijk} {z \choose k},
\]
with $c_{ijk} \in \ZZ_p$.
We now find that
\begin{eqnarray}
f_i(z) &=& s_i + \sum_{j=1}^{\infty} p^j 
\left( \sum_{k=0}^{N_j} c_{ijk} {z \choose k}\right) \nonumber\\
&=& s_0 + \sum_{k=0}^{\infty} b_{ik}  {z \choose k} \label{eq475}
\end{eqnarray}
in which
\[
b_{ik} := \sum_{j=1}^{\infty} p^j c_{ijk}
\]
is absolutely convergent $p$-adically, since each $c_{ijk} \in \ZZ_p$.
To show the series (\ref{eq475}) converges to an analytic map on
  $z \in \ZZ_p$, we must establish that $|b_{ik}|_{p}/|k!|_p \to 0$
as $k \to \infty$, i.e. that for any $j >0$ one has $b_{ik}/k! \in p^j \bZ$
for all sufficiently large $k$ \cite[Theorem 4.7, pp. 354]{Robert}.
To do this, we note that $c_{ijk}=0$ if $j<(k+1)/2$, hence
\[ b_{ik} := \sum_{j\ge (k+1)/2} p^j c_{ijk}.\]
It follows that $|b_{ik}|_p<p^{-(k+1)/2}$. Since $1/|k!|_p<p^{k/(p-1)}$,
we see 
that for all $p \ge 2$, there holds
$$
|b_{ik}/k!|_p \le p^{-k(\frac{1}{2}- \frac{1}{p-1}) - \frac{1}{2}} .
$$
Now if $p \ge 5$ this estimate gives 
\begin{equation}~\label{cbound}
|b_{ik}/k!|_p   \to 0, ~~\mbox{as}~~k \to \infty. 
\end{equation}
   Hence if $p \ge 5$ then $f_1(z) ,\ldots ,f_n(z) $ 
   each define analytic maps $f_i: \bZ \to \bZ$.\\

The argument above also establishes that
\[ 
f_i(z) \equiv g_{i,j}(z)~(\bmod~p^{j} \sRS).
\]
It then follows from property (3) above that
\[
f_i(z+1) \equiv H_i( f_1(z), ..., f_d(z)) (\bmod~ p^j \sRS)
\]
Since this holds for all $j \ge 1$, we conclude
that
\[
f_i(z+1) = H_i( f_1(z), ..., f_n(z))
\]
holds as p-adic analytic functions on $\ZZ_p$. This establishes (1).
Finally, we have
\[ 
f_{i}(0) = s_i + \sum_{j=1}^{\infty} p^j h_{i,j}(0) = s_i,
\]
which establishes (2). \qed \\
%
%

%
%

 We note that the conclusion of  Theorem~\ref{lem: 2nd} does not
always hold if $p=2$.  
\begin{exam}
Let $H(x)=-x$.  Then the following hold.
\begin{enumerate}
\item $H(1)=-1\equiv 1~(\bmod~ 2)$;
\item $\frac{dH}{dx} = -1 \equiv 1~(\bmod~ 2)$;
\item There does not exist a $2$-adic analytic function $f(z)$ satisfying $f(0)=1$ and $f(z+1)=H(f(z))$.
\end{enumerate}
\end{exam}
\noindent 
{\bf Proof.} It is straightforward to verify (1) and (2).  To show (3), we argue
by contradiction. So suppose that we have a $2$-adic analytic function $f(z)$ satisfying 
$f(0)=1$ and $f(z+1)=-f(z)$.  Then $f(2k)=1$ for every natural number $k$.  
By Strassman's theorem \cite{St} we see that $f(z)=1$ for every $n$.  But $f(1)=-1\not =1$, a contradiction.  Hence  (3) follows. \qed

We do not know of an example for which the conclusion of the statement of Theorem~\ref{lem: 2nd} fails when $p=3$.
\vskip 2mm

\noindent {\bf Proof of Theorem \ref{thm: introA}.}
By Lemma \ref{lem: p} and the remarks which follow it,
 we may assume that $\sigma\in {\rm Aut}(
\bA_{\bZ}^n)$, ${\bf q}=(q_1,\ldots ,q_n)\in \bA_{\bZ}^n$ 
and that $X$ is 
a subvariety of $\bA_{\bZ}^n$.
Let $\tau = \sigma^j$, where $j$ is as in Lemma \ref{lem: j}.  Let $i$ be 
between $0$ and $j-1$.
Either $\tau^m(\sigma^i({\bf q}))\in X$ for infinitely
many $m$ or there are only a finite number of $m$ for which this occurs.  
If the number of such $m$ is infinite, we claim that $\tau^m(\sigma^i({\bf q}))$
must be in $X$ for every $m$.
To see this we replace ${\bf q}$ with $\sigma^i({\bf q})$ so that
$\tau^n({\bf q})\in X$ for infinitely many $n$.
Write $\tau =(H_1,\ldots ,H_n)$, with 
$H_1,\ldots ,H_n$ all in ${\mathbb Z}_p[x_1,\ldots,H_n]$.
By construction, $H_i(q_1,\ldots ,q_n)\equiv q_i~(\bmod ~ p\bZ)$.
Hence by Theorem \ref{lem: 2nd}
there exist power series $f_1,\ldots ,f_n\in {\mathbb Z}_p[[z]]$
which converge in the closed unit ball of ${\mathbb Q}_p$ and 
which satisfy $f_i(0)=q_i$ and
$$f_i(z+1) \ = \ 
H_i(f_1(z),\ldots ,f_n(z))\qquad {\rm for}~0\le i\le n.$$  By construction,
$\tau^k({\bf q})=(f_1(k),\ldots ,f_n(k))$ for all $k\in {\mathbb Z}$.
By assumption the orbit of ${\bf q}$ under $\tau$ intersects $X$ infinitely often.
Write $X=Z(Q_1,\ldots ,Q_d)$.
For $1\le i\le d$, define
$$P_i(z)\ := \ Q_i(f_1(z),\ldots ,f_n(z))\in {\mathbb Q}_p[[z]].$$
Since $f_1,\ldots ,f_n$ converge on the closed unit disc and $Q_1,\ldots ,Q_d$ are
polynomials, $P_i(z)$ is
a $p$-adic power series which converges on the closed unit disc for $1\le i\le d$.
Moreover,
$P_i(m)=0$ whenever $\tau^{m}({\bf q})
\in X\subseteq Z(Q)$.  
Hence $P_i(z)$ has infinitely many zeros inside the closed unit ball.  
By Strassman's theorem \cite{St} $P_i(z)$ is identically zero.  
Thus $\tau^m({\bf q})$ is in $Z(Q_1,\ldots ,Q_d)=X$
for all integers $m$.  It follows that
$$\{ \sigma^{i+jm}({\bf q})~|~m\in {\mathbb Z}\} \ \subseteq \ 
X.$$
Hence for $0\le i\le j$, we have that 
$$\{m\in {\mathbb Z}~|~\sigma^{i+jm}({\bf q})\in X\}$$ is either a finite set or the
set of all integers.  The result now follows. \qed
\section{Proof of the main result}
In this section we use a theorem of Srinivas to deduce 
Theorem \ref{thm: intro}
from Theorem \ref{thm: introA}.
\begin{thm} (Srinivas) Let $A$ be a finitely generated algebra over an infinite
field $K$.  Then there exists a natural number $n=n(A)$ such that if $N>n$ and
$$f:K[x_1,\ldots ,x_N]\rightarrow A\qquad {\rm and}\qquad
g:K[x_1,\ldots ,x_N]\rightarrow A$$ are two surjections of $K$ algebras,
then there is an elementary $K$-algebra automorphism
$\phi:K[x_1,\ldots ,x_N]\rightarrow K[x_1,\ldots ,x_N]$ such that
$f=g\circ\phi$. \label{thm: Srinivas}
\end{thm}
{\bf Proof.} See Srinivas \cite[Theorem 2, p. 126]{Sr}. \qed
\vskip 1mm

Using this theorem we immediately obtain our main result.
\vskip 1mm
\noindent {\bf Proof of Theorem \ref{thm: intro}.} Let $A$ denote the ring of
regular functions of $Y$ and let $n$ be as in the statement of Theorem \ref{thm: Srinivas}. 
Since $A$ is finitely generated over $K$, we can find
 $N>n$ such that there exists a surjection
$$g:K[x_1,\ldots ,x_N]\rightarrow A.$$
The automorphism $\sigma$ of $Y$ induces a $K$-algebra automorphism $\psi: A\rightarrow
A$.  Let $f$ denote  
the composition 
$$K[x_1,\ldots ,x_N]\stackrel{g}{\rightarrow} A\stackrel{\psi}{\rightarrow} A.$$
By Theorem \ref{thm: Srinivas} there is an automorphism $\phi$ of
$K[x_1,\ldots ,x_N]$ such that $f=g\circ \phi$.  Equivalently,
$\psi\circ g=g\circ \phi$.
We have that $g$ induces an embedding $i:Y\hookrightarrow \bA^n_K$ and
$\phi$ induces
an automorphism $\tau$ of $\bA^n_K$; moreover, by construction, we have
$$\tau\circ i \ = \ i\circ \sigma.$$
Hence $\tau\big|_Y=\sigma$ and so  
$$\mathcal{S} \:= \ 
\{m~|~\sigma^m({\bf q})\in X\} \ = \ \{m~|~\tau^m({\bf q})\in X\}.$$
We may regard the subvariety $X$ of $Y$ as
a subvariety of $\bA^n_K$.  Then by Theorem \ref{thm: introA}, 
$\mathcal{S}$ is the union of
a finite set and the finite
union of complete doubly infinite arithmetic progressions. \qed     

\begin{cor} \label{cor: rog}
 Let $Y$ be an affine variety over a field $K$ of characteristic $0$,
let $\sigma$ be an automorphism
of $Y$ and let ${\bf q}$ be a point in $Y$.  Then the orbit of ${\bf q}$ under
$\sigma$ is Zariski dense if and only if for each proper subvariety $X$ of $Y$ the set
$\{m~|~\sigma^m({\bf q})\in X\}$ is a finite.
\end{cor}
{\bf Proof.}  Suppose first that $\{\sigma^m({\bf q})~|~m\in \mathbb{Z}\}$ 
is not dense in
$Y$.  Letting $X$ denote the Zariski closure of this set, we get an 
immediate contradiction since $\{m~|~\sigma^m({\bf q})\in X\}$ is a finite set
by assumption.
Suppose next
that there is some proper subvariety $X$ such that
$\{m~|~\sigma^m({\bf q})\in X\}$ is infinite.
Then by Theorem \ref{thm: intro}, 
there exist $i,j$ with $j\ge 1$
 such that $\sigma^{i+jm}({\bf q})\in X$ for all integers
$m$.
It follows that $$\{\sigma^m({\bf q})~|~m\in \mathbb{Z}\}\subseteq
\bigcup_{i=0}^{j-1} \sigma^i(X).$$
Thus the orbit is contained in a proper closed subset of $Y$ and so is not
Zariski dense.  \qed
\section{Examples}
The following example uses a \emph{Russell-Koras threefold} \cite{KR} to 
show that the iterates of a point with an infinite orbit
 under a non-linearizable
automorphism can intersect a non-linear subvariety $X$ of $\bA_{\mathbb C}^n$ in
a non-trivial arithmetic progression.  Furthermore, $X$ 
is not isomorphic to $\bA_{\mathbb C}^d$ for any $d\le n$.
\begin{exam}
Let
$\sigma\in {\rm Aut}(\bA_{\mathbb C}^4)$ be given by $$\sigma(a,b,c,d) \ = \ (a,b-3d^2-3da^2-a^4,
\omega c,
d+a^2)$$
where $\omega=e^{\pi i/k}$ (note that this automorphism depends on choice
of $k$). Let $X=Z(x+x^2y+z^2+t^3)=\{(a,b,c,d)\in \mathbb{C}^4~|~
a+a^2b+c^2+d^3=0\}\subseteq \bA_{\mathbb C}^4$ and let ${\bf q}=(-1,-1,1,1)$.
Then we have:
\begin{enumerate}
\item{
$X$ is not isomorphic to $\bA_{\mathbb C}^3$;}\item{ $\sigma$ is
not linearizable; that is, there is no automorphism $\rho$ of $\bA_{\mathbb C}^4$ such that
$\rho\circ \sigma \circ \rho^{-1}$ is an affine linear map;}\item{ $\{m~|~\sigma^m({\bf q})\in X\}$ is the
arithmetic progression $k\mathbb{Z}$} \item{${\bf q}$ has an infinite orbit under iteration of
$\sigma$.}\end{enumerate}
\end{exam}
{\bf Proof.} 
The fact that the affine subvariety $X\subseteq \bA^4_{\mathbb C}$ is not isomorphic to $\bA_{\mathbb C}^3$ 
is due to Makar-Limanov \cite{ML}.
Let $\tau(a,b,c,d)=(a,b-3d^2-3da^2-a^4,c,d+a^2)$.  
Then $\sigma$ and $\tau$ are 
automorphisms of $\bA_{\mathbb C}^4$.  
It is easy to check that $\tau$ restricts to an automorphism of
$X$.  Moreover, $\tau^{2k}=\sigma^{2k}$ and hence $\sigma^{2km}({\bf q})\in X$.  
For $j$ an odd multiple of $k$ we have that $\tau^j$ and $\sigma^j$ are the
same except that their third coordinates have opposite signs.  Hence
$\sigma^{km}({\bf q})\in X$.
For 
$j$ not a multiple of $k$, we have $\sigma^j({\bf q})$ differs from
$\tau^j({\bf q})$ in the $3$rd coordinate by a factor of $\omega^j$ and is
the same on the other coordinates.  Since 
$\omega^j\not =\pm 1$ we see that $\sigma^j({\bf q})\not \in X$ for $j$ not
a multiple of $k$.  
Notice that the fourth coordinate of $\sigma^j({\bf q})=1+j$ for all $j\in \mathbb{Z}$.  Hence
${\bf q}$ has an infinite orbit under iteration of $\sigma$.
Finally, to show that $\sigma$ is not linearizable, it is
sufficient to show that the automorphism
$\mu$ of $\bA_{\mathbb C}^3$ given by
$\mu(a,b,d)=(a,b-3d^2-3da^2-a^4,d+a^2)$ is not linearizable.  If $\mu$ were 
linearizable it would be of the form $\rho \circ L \circ \rho^{-1}$ where
$L$ is some affine linear automorphism and $\rho$ is some automorphism of 
$\bA^4_{\mathbb C}$.  Without loss of generality, we may assume
that $L$ is a linear map and that $\rho$ takes the origin to itself
since $\mu$ fixes the origin.  
But notice that the Jacobian of $\mu$ at the origin must be conjugate
to the Jacobian of $L$ at the origin.  But the Jacobian of $\mu$ at the origin
is the identity and the Jacobian of the linear map $L$ is $L$.  Hence $L$ is
conjugate to the identity and hence must be the identity map.  But then
$\mu$ must be the identity map, which is not the case. \qed
\vskip 1mm
 We now present an example showing that the conclusion to Theorem 
\ref{thm: intro} does not always hold in
positive characteristic.
This example is due to Lech 
\cite{Lech}, \cite{Bez2}.
\begin{prop} \label{exam: hoch}
Let $p$ be a prime and let $K=\overline{\mathbb{F}}_p(T)$.
Define $\sigma:\bA_K^2\rightarrow \bA_K^2$ by
$\sigma(x_1,x_2)=(Tx_1,(T+1)x_2)$.  Then for ${\bf q}=(1,1)$ and
$V=\{(a,b)\in K^2~|~x_2-x_1=1\}$, the forward orbit
$\sigma^j({\bf q})$ intersects $V$ infinitely often, but does not contain any
infinite arithmetic progression.
\end{prop}
{\bf Proof.}
Here $\sigma^j(1,1)=(T^j,(T+1)^j)$.  Thus
$\sigma^{p^d}(1,1)=(T^{p^d},T^{p^d}+1)\in V$.  
Notice, however, that if
$\sigma^{i+jm}(1,1)\subseteq Z(x_2-x_1-1)$ for some $i,j$ with $j\ge 1$ and
all nonnegative integers $m$, then
$(T+1)^{i+mj}-T^{i+mj}=1$ for all $m\ge 0$.
But since $T$ is transcendental over ${\mathbb F}_p$, the only way
for this to occur would be if $i+mj$ is a power of $p$ for all $m\ge 0$, which
is impossible.   The result follows. \qed
\vskip 2mm
It is natural to ask if there is a nice description of sets of the form
$$\{m\in \mathbb{Z}~|~\sigma^m({\bf q})\in X\}$$ for $X$ a subvariety of
$\mathbb{A}_K^n$, $\sigma\in {\rm Aut}(\bA_{K}^n)$, and ${\bf q}\in \bA_K^n$ when
we work over a field $K$ of characteristic $p>0$.
Examples seem to indicate that in characteristic $p>0$, these sets should
be $p$-\emph{automatic} sets \cite{AS}.  In fact, Derksen \cite{Derk} has 
recently given a positive characteristic analogue
 of the classical Skolem-Mahler-Lech theorem
 which shows that the zero set of a linear recurrence over a field
of positive characteristic $p$ is $p$-automatic; more specifically, he shows
that the zero set can be expressed as the union of a finite set, finitely
many arithmetic progressions, and a finite number of sets of the form
$$\{ c_1 p^{ak_1}+\cdots + c_d p^{ak_d} \} \cap \mathbb{N},$$ where
$a$ and $d$ are positive integers and $c_1,\ldots ,c_d$ are rational numbers.
In particular, Derksen's result implies
a result of B\'ezivin \cite{Bez2} which shows that in characteristic
$p>0$ the zero set of a linear recurrence either has
Banach density equal to $0$ or it contains an arithmetic progression.  
It is reasonable to conjecture that an analogue of Derksen's theorem holds
for affine varieties over fields of positive characteristic.
Derksen \cite{Derk} has constructed
pathological
examples which show that the zero set of a linear recurrence over a field of positive characteristic
 can be very complicated.  The following example, due to Derksen, shows the type of pathologies that can occur in 
positive characteristic.
\begin{exam} Let $p$ be prime and let $K=\mathbb{F}_p(u,t,s)$.
Define $\sigma:\bA_K^7\rightarrow \bA_K^7$ by
$$\sigma(x_1,x_2,\ldots ,x_7)=((u+t+s)x_1,(u+t)x_2,(u+s)x_3,(t+s)x_4,ux_5,tx_6,
sx_7).$$
  Then for ${\bf q}=(1,1,\ldots ,1)$ and
$$V=\{(a_1,\ldots ,a_7)\in K^7~|~a_1-a_2-a_3-a_4+a_5+a_6+a_7=0\},$$
 the set of natural numbers
$m$ for which 
$\sigma^m({\bf q})$ lies in $V$ is the set
$$\mathcal{S}
\ := \ \{p^i+p^j~|~i,j\in \mathbb{N}\}\cup \{p^i~|~i\in \mathbb{N}\}.$$
\end{exam}
{\bf Proof.}  The set of $m$ for which $\sigma^m({\bf q})\in V$ is just the
set of $m$ such that
$$(u+t+s)^m-(u+t)^m-(u+s)^m-(t+s)^m+u^m+t^m+s^m \ = \ 0.$$
First observe that if $m\not \in \mathcal{S}$, then the sum of the digits
in the base $p$ expansion of $m$ is at least $3$, and hence we can write
$m=p^i+p^j+p^kb$, with $i\le j\le k$ and $b$ not a multiple of $p$.
Then \begin{eqnarray*}
(u+t+s)^m & = & (u+t+s)^{p^i}(u+t+s)^{p^j}(u+t+s)^{p^kb}
\\ &= & (u^{p^i}+t^{p^i}+s^{p^i})(u^{p^j}+t^{p^j}+s^{p^j})(u^{p^k}+t^{p^k}+s^{p^k})^b, 
\end{eqnarray*}
which has nonzero coefficient of $u^{p^i}t^{p^j}s^{p^kb}$ since
$u,t$ and $s$ are algebraically independent.  Consequently,
$$(u+t+s)^m-(u+t)^m-(u+s)^m-(t+s)^m+u^m+t^m+s^m \ \not = \ 0$$ if
$m\not \in \mathcal{S}$.

On the other hand, if $m=p^i$, then we have
\begin{eqnarray*}
&~& 
(u+t+s)^m-(u+t)^m-(u+s)^m-(t+s)^m+u^m+t^m+s^m \\
& = & (u^{p^i}+t^{p^i}+s^{p^i})
- (u^{p^i}+t^{p^i}) - (u^{p^i}+t^{p^i}) - (t^{p^i}+s^{p^i}) + u^{p^i}+t^{p^i}+s^{p^i} \\
& = & 0.\end{eqnarray*}
Similarly, if $m=p^i+p^j$, an easy computation shows that
$$(u+t+s)^m-(u+t)^m-(u+s)^m-(t+s)^m+u^m+t^m+s^m \ = \ 0.$$
  \qed
\section{Concluding remarks}
Throughout this paper, we have been dealing with automorphisms of an affine
variety $Y$.  There is no reason why one could not instead study endomorphisms
of affine varieties.  We raise the following question. 
\begin{quest} Let $Y$ be an
affine variety over a field $K$ of characteristic $0$ 
with a subvariety $X$ and suppose that $\tau:Y\rightarrow Y$ and ${\bf q}\in Y$.
Is it true that the set 
of $m\ge 0$ such that $\tau^m({\bf q})\in X$ is the union of a finite set and
a finite number of complete arithmetic progressions?
\label{conj: 1}\end{quest}
The methods of this paper do not directly apply to this question.  Our results give some weak implications
for the Jacobian conjecture.  Lemmas \ref{lem: j} and \ref{lem: comp} avoid the hypothesis
that $\sigma$ is an automorphism, instead using hypotheses about it being 
surjective and having constant Jacobian determinant.  
With these Lemmas one can deduce, for instance, that if ${\bf F}:\bZ^n\rightarrow
\bZ^n$ is a surjective polynomial map with constant Jacobian determinant, then
${\bf F}$ has a polynomial inverse.   

Denis \cite{Denis} studied the orbit of a point in projective space under a morphism, showing
that in characteristic $0$ if the orbit has infinite intersection with 
a subvariety that is minimal with respect to this property, then 
the subvariety is stable under some power of the morphism.  
This result holds in our setting as well.
Suppose that
 $Y$ is an affine variety over a field of characteristic $0$ 
and suppose that $\sigma$ is an automorphism
of $Y$.  If ${\bf q}$ is a point that has an infinite orbit under $Y$, then
the Zariski closure of this orbit is some subvariety $X$.  The orbit must 
intersect some irreducible subvariety of $X$ in an infinite set.  We then
see that some power of $\sigma$ acts as an automorphism of this irreducible
subvariety.

Finally, we note that there are examples from complex dynamics of automorphisms of
the affine plane $\bA^2_{\mathbb C}$ with
strange ``attractor sets.''  The most famous examples are the \emph{H\'enon mappings} \cite{Hen}, \cite{For}
given by
$${\bf H}:(x,y)\rightarrow (p(x)-ay,x),$$ 
where $a$ is a real number and $p(x)$ is a polynomial of degree at least two  with real coefficients.
For these maps there exists an associated ``Julia set'' $J$ in $\bA^2_{\mathbb C}$
invariant under ${\bf H}$ (see Bedford and Smillie \cite{BS1}).
Theorem \ref{thm: intro} may say something about how $J$ intersects a curve in $\bA^2_{\mathbb C}$.
\section{Acknowledgments}
I thank Jeff Lagarias for many invaluable comments.
I also thank Toby Stafford, Dan Rogalski, Harm Derksen, 
and Peter Russell for many
helpful discussions.  Finally, I thank the referee for many helpful comments.

\end{document}